# Definitions Derived from Neutrosophics

Florentin Smarandache
Department of Mathematics
University of New Mexico
200 College Road
Gallup, NM 87301, USA

*Abstract*: Thirty-three new definitions are presented, derived from neutrosophic set, neutrosophic probability, neutrosophic statistics, and neutrosophic logic.
Each one is independent, short, with references and cross references like in a dictionary style.

*Introduction*:
As an addenda to [1], [3], [5-7] we display the below unusual extension of definitions resulted from neutrosophics in the Set Theory, Probability, and Logic. Some of them are listed in the Dictionary of Computing [2]. Further development of these definitions (including properties, applications, etc.) is in our research plan.

*1. Definitions of New Sets*

=====================================================

## 1.1. Neutrosophic Set:

<logic, mathematics> A set which generalizes many existing classes of sets, especially the fuzzy set.

Let U be a universe of discourse, and M a set included in U. An element x from U is noted, with respect to the set M, as x(T,I,F), and belongs to M in the following way: it is T% in the set (membership appurtenance), I% indeterminate (unknown if it is in the set), and F% not in the set (non-membership); here T,I,F are real standard or non-standard subsets, included in the non-standard unit interval $]^-0, 1^+[$, representing truth, indeterminacy, and falsity percentages respectively.

Therefore: $^-0 \le \inf(T) + \inf(I) + \inf(F) \le \sup(T) + \sup(I) + \sup(F) \le 3^+$.

Generalization of {classical set}, {fuzzy set}, {intuitionistic set}, {paraconsistent set}, {faillibilist set}, {paradoxist set}, {tautological set}, {nihilist set}, {dialetheist set}, {trivialist}.

Related to {neutrosophic logic}.

{ ref. Florentin Smarandache, "A Unifying Field in Logics. Neutrosophy: Neutrosophic Probability, Set, and Logic", American Research Press, Rehoboth, 1999; (http://www.gallup.unm.edu/~smarandache/NeutrosophicSet.pdf, http://www.gallup.unm.edu/~smarandache/FirstNeutConf.htm, http://www.gallup.unm.edu/~smarandache/neut-ad.htm) }

=======================================================

## 1.2. Intuitionistic Set:

<logic, mathematics> A set which provides incomplete information on its elements.

A class of {neutrosophic set} in which every element x is incompletely known, i.e. x(T,I,F) such that sup(T)+sup(I)+sup(F)<1; here T,I,F are real standard or non-standard subsets, included in the non-standard unit interval $]^-0, 1^+[$, representing truth, indeterminacy, and falsity percentages respectively.

Contrast with {paraconsistent set}.

Related to {intuitionistic logic}.

{ ref. Florentin Smarandache, "A Unifying Field in Logics.
Neutrosophy: Neutrosophic Probability, Set, and Logic",
American Research Press, Rehoboth, 1999;
(http://www.gallup.unm.edu/~smarandache/FirstNeutConf.htm,
 http://www.gallup.unm.edu/~smarandache/neut-ad.htm) }

=====================================================

1.3. Paraconsistent Set:

<logic, mathematics> A set which provides paraconsistent information on its elements.

A class of {neutrosophic set} in which every element x(T,I,F) has the property that sup(T)+sup(I)+sup(F)>1;
here T,I,F are real standard or non-standard subsets, included in the non-standard unit interval $]^-0, 1^+[$, representing truth, indeterminacy, and falsity percentages respectively.

Contrast with {intuitionistic set}.

Related to {paraconsistent logic}.

{ ref. Florentin Smarandache, "A Unifying Field in Logics.
Neutrosophy: Neutrosophic Probability, Set, and Logic",
American Research Press, Rehoboth, 1999;
(http://www.gallup.unm.edu/~smarandache/FirstNeutConf.htm,
 http://www.gallup.unm.edu/~smarandache/neut-ad.htm) }

=====================================================

1.4. Faillibilist Set:

<logic, mathematics> A set whose elements are uncertain.

A class of {neutrosophic set} in which every element x has a percentage of indeterminacy, i.e. x(T,I,F) such that inf(I)>0;
here T,I,F are real standard or non-standard subsets, included in the non-standard unit interval $]^-0, 1^+[$, representing truth, indeterminacy, and falsity percentages respectively.

Related to {faillibilism}.

{ ref. Florentin Smarandache, "A Unifying Field in Logics.
Neutrosophy: Neutrosophic Probability, Set, and Logic",

American Research Press, Rehoboth, 1999;
(http://www.gallup.unm.edu/~smarandache/FirstNeutConf.htm,
 http://www.gallup.unm.edu/~smarandache/neut-ad.htm) }

==================================================

1.5. Paradoxist Set:

<logic, mathematics> A set which contains and doesn't contain
itself at the same time.

A class of {neutrosophic set} in which every element x(T,I,F) has
the form x(1,I,1), i.e. belongs 100% to the set and doesn't
belong 100% to the set simultaneously;
here T,I,F are real standard or non-standard subsets, included in
the non-standard unit interval $]^-0, 1^+[$, representing truth,
indeterminacy, and falsity percentages respectively.

Related to {paradoxism}.

{ ref. Florentin Smarandache, "A Unifying Field in Logics.
Neutrosophy: Neutrosophic Probability, Set, and Logic",
American Research Press, Rehoboth, 1999;
(http://www.gallup.unm.edu/~smarandache/FirstNeutConf.htm,
 http://www.gallup.unm.edu/~smarandache/neut-ad.htm) }

==================================================

1.6. Pseudo-Paradoxist Set:

<logic, mathematics> A set which totally contains and partially doesn't contain
itself at the same time,
or partially contains and totally doesn't contain itself at the same time.

A class of {neutrosophic set} in which every element x(T,I,F) has
the form x(1,I,F) with 0<inf(F)   sup(F)<1 or x(T,I,1) with 0<inf(T)   sup(T)<1,
i.e. belongs 100% to the set and doesn't belong F% to the set simultaneously, with
0<inf(F)   sup(F)<1,
or belongs T% to the set and doesn't belong 100% to the set simultaneously, with
0<inf(T)   sup(T)<1;
here T,I,F are real standard or non-standard subsets, included in
the non-standard unit interval $]^-0, 1^+[$, representing truth,
indeterminacy, and falsity percentages respectively.

Related to {pseudo-paradoxism}.

{ ref. Florentin Smarandache, "A Unifying Field in Logics.
Neutrosophy: Neutrosophic Probability, Set, and Logic",
American Research Press, Rehoboth, 1999;
(http://www.gallup.unm.edu/~smarandache/FirstNeutConf.htm,
 http://www.gallup.unm.edu/~smarandache/neut-ad.htm) }

===================================================

1.7. Tautological Set:

<logic, mathematics> A set whose elements are absolutely
determined in all possible worlds.

A class of {neutrosophic set} in which every element x has the
form $x(1^+, {}^-0, {}^-0)$, i.e. absolutely belongs to the set;
here T,I,F are real standard or non-standard subsets, included
in the non-standard unit interval $]^-0, 1^+[$, representing truth,
indeterminacy, and falsity percentages respectively.

Contrast with {nihilist set} and {nihilism}.

Related to {tautologism}.

{ ref. Florentin Smarandache, "A Unifying Field in Logics.
Neutrosophy: Neutrosophic Probability, Set, and Logic",
American Research Press, Rehoboth, 1999;
(http://www.gallup.unm.edu/~smarandache/FirstNeutConf.htm,
 http://www.gallup.unm.edu/~smarandache/neut-ad.htm) }

===================================================

1.8. Nihilist Set:

<logic, mathematics> A set whose elements absolutely
don't belong to the set in all possible worlds.

A class of {neutrosophic set} in which every element x has the
form $x({}^-0, {}^-0, 1^+)$, i.e. absolutely doesn't belongs to the set;
here T,I,F are real standard or non-standard subsets, included
in the non-standard unit interval $]^-0, 1^+[$, representing truth,
indeterminacy, and falsity percentages respectively.

The empty set is a particular set of {nihilist set}.

Contrast with {tautological set}.

Related to {nihilism}.

{ ref. Florentin Smarandache, "A Unifying Field in Logics.
Neutrosophy: Neutrosophic Probability, Set, and Logic",
American Research Press, Rehoboth, 1999;
(http://www.gallup.unm.edu/~smarandache/FirstNeutConf.htm,
 http://www.gallup.unm.edu/~smarandache/neut-ad.htm) }

========================================================

1.9. Dialetheist Set:

 <logic, mathematics> /di:-al-u-theist/ A set which contains at
least one element which also belongs to its complement.

 A class of {neutrosophic set} which models a situation
 where the intersection of some disjoint sets is not empty.

 There is at least one element x(T,I,F) of the dialetheist set
 M which belongs at the same time to M and to the set C(M),
 which is the complement of M;
 here T,I,F are real standard or non-standard subsets, included in the
 non-standard unit interval $]^-0, 1^+[$, representing truth,
 indeterminacy, and falsity percentages respectively.

 Contrast with {trivialist set}.

 Related to {dialetheism}.

 { ref. Florentin Smarandache, "A Unifying Field in Logics.
 Neutrosophy: Neutrosophic Probability, Set, and Logic",
 American Research Press, Rehoboth, 1999;
 (http://www.gallup.unm.edu/~smarandache/FirstNeutConf.htm,
  http://www.gallup.unm.edu/~smarandache/neut-ad.htm) }

========================================================

1.10. Trivialist Set:

 <logic, mathematics>  A set all of whose  elements also belong
to its complement.

 A class of {neutrosophic set} which models a situation
 where the intersection of any disjoint sets is not empty.

 Every element x(T,I,F) of the trivialist set M belongs at the

same time to M and to the set C(M), which is the
complement of M;
here T,I,F are real standard or non-standard subsets,
included in the non-standard unit interval $]^-0, 1^+[$, representing
truth, indeterminacy, and falsity percentages respectively.

Contrast with {dialetheist set}.

Related to {trivialism}.

{ ref. Florentin Smarandache, "A Unifying Field in Logics.
Neutrosophy: Neutrosophic Probability, Set, and Logic",
American Research Press, Rehoboth, 1999;
(http://www.gallup.unm.edu/~smarandache/FirstNeutConf.htm,
 http://www.gallup.unm.edu/~smarandache/neut-ad.htm) }

========================================================

2. Definitions of New Probabilities and Statistics

========================================================

2.1. Neutrosophic Probability:

<probability> The probability that an event occurs is (T, I, F),
where T,I,F are real standard or non-standard subsets, included in the
non-standard unit interval $]^-0, 1^+[$, representing truth,
indeterminacy, and falsity percentages respectively.

Therefore: $^-0 \le \inf(T) + \inf(I) + \inf(F) \le \sup(T) + \sup(I) + \sup(F) \le 3^+$.

Generalization of {classical probability} and {imprecise probability},
{intuitionistic probability}, {paraconsistent probability}, {faillibilist
probability}, {paradoxist probability}, {tautological probability},
{nihilistic probability}, {dialetheist probability}, {trivialist probability}.

Related with {neutrosophic set} and {neutrosophic logic}.

The analysis of neutrosophic events is called Neutrosophic Statistics.

{ ref. Florentin Smarandache, "A Unifying Field in Logics.
Neutrosophy: Neutrosophic Probability, Set, and Logic",
American Research Press, Rehoboth, 1999;
(http://www.gallup.unm.edu/~smarandache/FirstNeutConf.htm,
 http://www.gallup.unm.edu/~smarandache/neut-ad.htm) }

=====================================================

2.2. Intuitionistic Probability:

<probability> The probability that an event occurs is (T, I, F),
where T,I,F are real standard or non-standard subsets, included in the
non-standard unit interval $]^{-}0, 1^{+}[$, representing truth,
indeterminacy, and falsity percentages respectively,
and $n\_sup = sup(T)+sup(I)+sup(F) < 1$,
i.e. the probability is incompletely calculated.

Contrast with {paraconsistent probability}.

Related to {intuitionistic set} and {intuitionistic logic}.

The analysis of intuitionistic events is called Intuitionistic Statistics.

{ ref. Florentin Smarandache, "A Unifying Field in Logics.
 Neutrosophy: Neutrosophic Probability, Set, and Logic",
 American Research Press, Rehoboth, 1999;
(http://www.gallup.unm.edu/~smarandache/FirstNeutConf.htm,
 http://www.gallup.unm.edu/~smarandache/neut-ad.htm) }

=====================================================

2.3. Paraconsistent Probability:

<probability> The probability that an event occurs is (T, I, F),
where T,I,F are real standard or non-standard subsets, included in the
non-standard unit interval $]^{-}0, 1^{+}[$, representing truth,
indeterminacy, and falsity percentages respectively,
and $n\_sup = sup(T)+sup(I)+sup(F) > 1$,
i.e. contradictory information from various sources.

Contrast with {intuitionistic probability}.

Related to {paraconsistent set} and {paraconsistent logic}.

The analysis of paraconsistent events is called
Paraconsistent Statistics.

{ ref. Florentin Smarandache, "A Unifying Field in Logics.
 Neutrosophy: Neutrosophic Probability, Set, and Logic",
 American Research Press, Rehoboth, 1999;
 (http://www.gallup.unm.edu/~smarandache/FirstNeutConf.htm,
  http://www.gallup.unm.edu/~smarandache/neut-ad.htm) }

=======================================================

2.4. Faillibilist Probability:

<probability> The probability that an event occurs is (T, I, F),
where T,I,F are real standard or non-standard subsets, included in the
non-standard unit interval $]^{-}0, 1^{+}[$, representing truth,
indeterminacy, and falsity percentages respectively,
and inf(I) > 0,
i.e. there is some percentage of indeterminacy in calculation.

Related to {faillibilist set} and {faillibilism}.

The analysis of faillibilist events is called Faillibilist Statistics.

{ ref. Florentin Smarandache, "A Unifying Field in Logics.
 Neutrosophy: Neutrosophic Probability, Set, and Logic",
 American Research Press, Rehoboth, 1999;
 (http://www.gallup.unm.edu/~smarandache/FirstNeutConf.htm,
  http://www.gallup.unm.edu/~smarandache/neut-ad.htm) }

=======================================================

2.5. Paradoxist Probability:

<probability> The probability that an event occurs is (1, I, 1),
where I is a standard or non-standard subset, included in the
non-standard unit interval $]^{-}0, 1^{+}[$, representing indeterminacy.

Paradoxist probability is used for paradoxal events (i.e. which
may occur and may not occur simultaneously).

Related to {paradoxist set} and {paradoxism}.

The analysis of paradoxist events is called Paradoxist Statistics.

{ ref. Florentin Smarandache, "A Unifying Field in Logics.
 Neutrosophy: Neutrosophic Probability, Set, and Logic",
 American Research Press, Rehoboth, 1999;
(http://www.gallup.unm.edu/~smarandache/FirstNeutConf.htm,
  http://www.gallup.unm.edu/~smarandache/neut-ad.htm) }

=======================================================

2.6. Pseudo-Paradoxist Probability:

<probability> The probability that an event occurs is either (1, I, F) with
0<inf(F) $\leq$ sup(F)<1, or (T, I, 1) with 0<inf(T) $\leq$ sup(T)<1,
where T,I,F are standard or non-standard subset, included in the
non-standard unit interval $]^-0, 1^+[$, representing the truth, indeterminacy, and
falsity percentages respectively.

Pseudo-Paradoxist probability is used for pseudo-paradoxal events (i.e. which
may certainly occur and may not partially occur simultaneously,
or may partially occur and may not certainly occur simultaneously).

Related to {pseudo-paradoxist set} and {pseudo-paradoxism}.

The analysis of pseudo-paradoxist events is called Pseudo-Paradoxist Statistics.

{ ref. Florentin Smarandache, "A Unifying Field in Logics.
 Neutrosophy: Neutrosophic Probability, Set, and Logic",
 American Research Press, Rehoboth, 1999;
(http://www.gallup.unm.edu/~smarandache/FirstNeutConf.htm,
  http://www.gallup.unm.edu/~smarandache/neut-ad.htm) }

====================================================

2.7. Tautological Probability:

<probability> The probability that an event occurs is more than one,
i.e. $(1^+, {}^-0, {}^-0)$.

Tautological probability is used for universally sure events (in all
possible worlds, i.e. do not depend on time, space, subjectivity, etc.).

Contrast with {nihilistic probability} and {nihilism}.

Related to {tautological set} and {tautologism}.

The analysis of tautological events is called Tautological Statistics.

{ ref. Florentin Smarandache, "A Unifying Field in Logics.
 Neutrosophy: Neutrosophic Probability, Set, and Logic",
 American Research Press, Rehoboth, 1999;
(http://www.gallup.unm.edu/~smarandache/FirstNeutConf.htm,
  http://www.gallup.unm.edu/~smarandache/neut-ad.htm) }

====================================================

2.8. Nihilist Probability:

<probability> The probability that an event occurs is less than zero, i.e. ($^-0$, $^-0$, $1^+$).

Nihilist probability is used for universally impossible events (in all possible worlds, i.e. do not depend on time, space, subjectivity, etc.).

Contrast with {tautological probability} and {tautologism}.

Related to {nihilist set} and {nihilism}.

The analysis of nihilist events is called Nihilist Statistics.

{ ref. Florentin Smarandache, "A Unifying Field in Logics. Neutrosophy: Neutrosophic Probability, Set, and Logic", American Research Press, Rehoboth, 1999; (http://www.gallup.unm.edu/~smarandache/FirstNeutConf.htm, http://www.gallup.unm.edu/~smarandache/neut-ad.htm) }

========================================================

2.9. Dialetheist Probability:

<probability> /di:-al-u-theist/ A probability space where at least one event and its complement are not disjoint.

A class of {neutrosophic probability} that models a situation where the intersection of some disjoint events is not empty.

Here, similarly, the probability of an event to occur is (T, I, F), where T,I,F are real standard or non-standard subsets, included in the non-standard unit interval $]^-0, 1^+[$, representing truth, indeterminacy, and falsity percentages respectively.

Contrast with {trivialist probability}.

Related to {dialetheist set} and {dialetheism}.

The analysis of dialetheist events is called Dialetheist Statistics.

{ ref. Florentin Smarandache, "A Unifying Field in Logics. Neutrosophy: Neutrosophic Probability, Set, and Logic", American Research Press, Rehoboth, 1999; (http://www.gallup.unm.edu/~smarandache/FirstNeutConf.htm, http://www.gallup.unm.edu/~smarandache/neut-ad.htm) }

==========================================================

2.10. Trivialist Probability:

<probability> A probability space where every event and its complement are not disjoint.

A class of {neutrosophic probability} which models a situation where the intersection of any disjoint events is not empty.

Here, similarly, the probability of an event to occur is (T, I, F), where T,I,F are real standard or non-standard subsets, included in the non-standard unit interval $]^-0, 1^+[$, representing truth, indeterminacy, and falsity percentages respectively.

Contrast with {dialetheist probability}.

Related to {trivialist set} and {trivialism}.

The analysis of trivialist events is called Trivialist Statistics.

{ ref. Florentin Smarandache, "A Unifying Field in Logics. Neutrosophy: Neutrosophic Probability, Set, and Logic", American Research Press, Rehoboth, 1999; (http://www.gallup.unm.edu/~smarandache/FirstNeutConf.htm, http://www.gallup.unm.edu/~smarandache/neut-ad.htm) }

==========================================================

*3. Definitions of New Logics*

==========================================================

3.1. Neutrosophic Logic:

<logic, mathematics> A logic which generalizes many existing classes of logics, especially the fuzzy logic.

In this logic each proposition is estimated to have the percentage of truth in a subset T, the percentage of indeterminacy in a subset I, and the percentage of falsity in a subset F;
here T,I,F are real standard or non-standard subsets, included in the non-standard unit interval $]^-0, 1^+[$, representing truth, indeterminacy, and falsity percentages respectively.

Therefore: $^-0 \leq \inf(T) + \inf(I) + \inf(F) \leq \sup(T) + \sup(I) + \sup(F) \leq 3^+$.

Generalization of {classical or Boolean logic}, {fuzzy logic},
{multiple-valued logic}, {intuitionistic logic}, {paraconsistent logic},
{faillibilist logic, or failibilism}, {paradoxist logic, or paradoxism},
{pseudo-paradoxist logic, or pseudo-paradoxism}, {tautological logic, or
tautologism}, {nihilist logic, or nihilism}, {dialetheist logic, or dialetheism},
{trivialist logic, or trivialism}.

 Related to {neutrosophic set}.

 { ref. Florentin Smarandache, "A Unifying Field in Logics.
 Neutrosophy: Neutrosophic Probability, Set, and Logic",
 American Research Press, Rehoboth, 1999;
(http://www.gallup.unm.edu/~smarandache/NeutrosophicLogic.pdf,
 http://www.gallup.unm.edu/~smarandache/FirstNeutConf.htm,
 http://www.gallup.unm.edu/~smarandache/neut-ad.htm) }

=====================================================

3.2. Paradoxist Logic (or Paradoxism):

 <logic, mathematics> A logic devoted to paradoxes, in which each
 proposition has the logical vector value (1, I, 1);
 here I is a real standard or non-standard subset, included in the
 non-standard unit interval $]^-0, 1^+[$, representing the indeterminacy.

As seen, each paradoxist (paradoxal) proposition is true and false
simultaneously.

 Related to {paradoxist set}.

 { ref. Florentin Smarandache, "A Unifying Field in Logics.
 Neutrosophy: Neutrosophic Probability, Set, and Logic",
 American Research Press, Rehoboth, 1999;
(http://www.gallup.unm.edu/~smarandache/NeutrosophicLogic.pdf,
 http://www.gallup.unm.edu/~smarandache/FirstNeutConf.htm,
 http://www.gallup.unm.edu/~smarandache/neut-ad.htm) }

   =====================================================

3.3. Pseudo-Paradoxist Logic (or Pseudo-Paradoxism):

 <logic, mathematics> A logic devoted to pseudo-paradoxes,
 in which each proposition has the logical vector value:
 either (1, I, F), with $0<\inf(F)\ \sup(F)<1$,
 or (T, I, 1), with $0<\inf(T)\ \sup(T)<1$;

here I is a real standard or non-standard subset, included in the
non-standard unit interval $]^-0, 1^+[$, representing the indeterminacy.

As seen, each pseudo-paradoxist (pseudo-paradoxal) proposition is:
either totally true and partially false simultaneously,
or partially true and totally false simultaneously.

Related to {pseudo-paradoxist set}.

 { ref. Florentin Smarandache, "A Unifying Field in Logics.
 Neutrosophy: Neutrosophic Probability, Set, and Logic",
 American Research Press, Rehoboth, 1999;
(http://www.gallup.unm.edu/~smarandache/NeutrosophicLogic.pdf,
 http://www.gallup.unm.edu/~smarandache/FirstNeutConf.htm,
 http://www.gallup.unm.edu/~smarandache/neut-ad.htm) }

   ========================================================

3.4. Tautological Logic (or Tautologism):

 <logic, mathematics> A logic devoted to tautologies, in which each
 proposition has the logical vector value $(1^+, {}^-0, {}^-0)$.

As seen, each tautological proposition is absolutely true (i. e, true in all
possible worlds).

Related to {tautological set}.

 { ref. Florentin Smarandache, "A Unifying Field in Logics.
 Neutrosophy: Neutrosophic Probability, Set, and Logic",
 American Research Press, Rehoboth, 1999;
(http://www.gallup.unm.edu/~smarandache/NeutrosophicLogic.pdf,
 http://www.gallup.unm.edu/~smarandache/FirstNeutConf.htm,
 http://www.gallup.unm.edu/~smarandache/neut-ad.htm) }

   ========================================================

General References: